\documentclass{lmcs} %%% last changed 2014-08-20
\pdfoutput=1

\usepackage{lastpage}
\lmcsdoi{15}{2}{9}
\lmcsheading{}{\pageref{LastPage}}{}{}%
{Oct.~27,~2017}{May~08,~2019}{}


\usepackage{hyperref}
\usepackage{amssymb}
\usepackage{amsmath}

\usepackage[all]{xy}
\usepackage{graphicx}
\theoremstyle{plain} %\crefname{satz}{Satz}{S\"atze}
\newtheorem{them}[thm] {Theorem}
\newtheorem{Exa} [thm] {Example}
\newtheorem{Rem} [thm] {Remark}
\newcommand{\Ana}{{\it\Sigma}^{1}_{1}}

\newcommand{\Ca}{{\it\Pi}^{1}_{1}}
\newcommand{\Boraone}{{\it\Sigma}^{0}_{1}}
\newcommand{\ca}{{\bf\Pi}^{1}_{1}}
\newcommand{\Borel}{{\it\Delta}^{1}_{1}}
\newcommand{\ana}{{\bf\Sigma}^{1}_{1}}

\newcommand{\boraone}{{\bf\Sigma}^{0}_{1}}
\newcommand{\boratwo}{{\bf\Sigma}^{0}_{2}}

\newcommand{\boraxi}{{\bf\Sigma}^{0}_{\xi}}

\newcommand{\borone}{{\bf\Delta}^{0}_{1}}
\newcommand{\bortwo}{{\bf\Delta}^{0}_{2}}

\newcommand{\bormone}{{\bf\Pi}^{0}_{1}}
\newcommand{\Bormone}{{\it\Pi}^{0}_{1}}

\newcommand{\bormtwo}{{\bf\Pi}^{0}_{2}}

\newcommand{\bormxi}{{\bf\Pi}^{0}_{\xi}}

\newcommand{\borxi}{{\bf\Delta}^{0}_{\xi}}
\newcommand{\borme}{{\bf\Pi}^{0}_{\eta}}

\newcommand{\om}{\omega}

\newcommand{\nl}{\newline}
\newcommand{\fa}{\forall}

\newcommand{\Si}{\Sigma}

\newcommand{\Sio}{\Sigma^\omega}
\newcommand{\ra}{\rightarrow}

\newcommand{\hs}{\hspace{12mm}

\noindent}

\begin{document}

\title{Polishness of some topologies related to word or tree automata}

\author[Olivier Carton]{Olivier Carton}	%required
\address{Universit\'e Paris Diderot, Institut de Recherche en Informatique Fondamentale, UMR 8243, Case 7014,  75 205 Paris Cedex 13, France}	%required
\email{Olivier.Carton@irif.fr}  %optional
%\thanks{thanks 1, optional.}	%optional

\author[Olivier Finkel]{Olivier Finkel}	%optional
\address{Equipe de Logique Math\'ematique, Institut de Math\'ematiques de Jussieu-Paris Rive Gauche, UMR7586, CNRS et Universit\'e Paris 7, France}	%optional
\email{finkel@math.univ-paris-diderot.fr}  %optional
%\thanks{thanks 2, optional.}	%optional

\author[Dominique Lecomte]{Dominique Lecomte}	%optional
\address{Projet Analyse Fonctionnelle, Institut de Math\'ematiques de Jussieu-Paris Rive Gauche, Universit\'e Paris 6, France, ~~ {\it and}    Universit\'e de Picardie, I.U.T. de l'Oise, site de Creil, France} %optional
\email{dominique.lecomte@upmc.fr}
%\urladdr{name3@url3\quad\rm{(optionally, a web-page can be specified)}}  %optional
%\thanks{thanks 3, optional.}	%optional

\begin{abstract}
  We prove that the B\"uchi topology and the automatic topology are Polish. We also show that this cannot be fully extended to the case of the space of infinite labelled binary trees; in particular the B\"uchi and the Muller topologies are not Polish in this case.  
\end{abstract}

\maketitle

\section{Introduction}\indent

This paper is a contribution to the study of the interactions between
descriptive set theory and theoretical computer science. These interactions
have already been the subject of many studies, see for instance
\cite{Landweber69,Wagner79,LescowThomas,Staiger97,Selivanov03b,NW03,PerrinPin,Serre04b,Fin-mscs06,Selivanov08,Fin-ICST,Fink-Sim,Murlak08,Fink-Lec2,DFR4,Skrzypczak16,CavallariMS17,MMS}.

In particular, the theory of automata reading infinite words, which is
closely related to infinite games, is now a rich theory which is used for
the specification and the verification of non-terminating systems, see
\cite{2001automata,PerrinPin}.  The space $\Sigma^{\mathbb{N}}$ of infinite
words over a finite alphabet $\Sigma$, equipped with the usual Cantor
topology $\tau_C$, is a natural place to study the topological complexity of the $\omega$-languages
accepted by various kinds of automata. In particular, it is interesting to locate them with respect to the Borel and
the projective hierarchies.

However, as noticed in \cite{SchwarzS10} by Schwarz and Staiger and in
\cite{HoffmannS15} by Hoffmann and Staiger, it turns out that for several
purposes some other topologies on the space $\Sigma^{\mathbb{N}}$ are useful,
for instance for studying fragments of the first-order logic over infinite
words, or for a topological characterization of the random infinite words (see also \cite{HoffmannSS16}).  In
particular, Schwarz and Staiger studied four topologies on the space
$\Sigma^{\mathbb{N}}$ of infinite words over a finite alphabet $\Sigma$, 
which are all related to automata, and refine the Cantor topology on
$\Sigma^{\mathbb{N}}$: the B\"uchi topology, the automatic topology, the
alphabetic topology, and the strong alphabetic topology.

Recall that a topological space is Polish if and only if it is separable,
i.e., contains a countable dense subset, and its topology is induced
by a complete metric.  Classical descriptive set theory is about the
topological complexity of the definable subsets of the Polish topological spaces,
as well as the study of some hierarchies of topological complexity (see
\cite{Kechris94,PerrinPin} for the basic notions). The analytic sets, which are
the projections of the Borel sets, are of particular importance. Similar
hierarchies of complexity are studied in effective descriptive set theory,
which is based on the theory of recursive functions (see
\cite{Moschovakis09} for the basic notions). The effective analytic subsets of
the Cantor space $(2^{\mathbb{N}},\tau_C)$ are highly related to
theoretical computer science, in the sense that they coincide with the sets
recognized by some special kind of Turing machine (see \cite{Staiger99}).
  
We now give details about some topologies that we investigate in this paper. Let $\Sigma$ be a finite alphabet with at least two symbols. We consider the following topologies on $\Sigma^{\mathbb{N}}$. 

\begin{itemize}
\item the \emph{B\"uchi topology} $\tau_B$, generated by the set
  $\mathbb{B}_B$ of $\omega$-regular languages,\smallskip
 
\item the \emph{automatic topology} $\tau_A$, generated by the set
  $\mathbb{B}_A$ of $\tau_C$-closed $\omega$-regular languages (this
  topology is remarkable because any $\tau_C$-closed (or even $\tau_C$-${\bf \Pi}^0_2$) $\omega$-regular language is accepted by some deterministic B\"uchi automaton, \cite{PerrinPin}),\smallskip
 
\item the \emph{topology} $\tau_\delta$, generated by the set $\mathbb{B}_\delta$ of languages accepted by some unambiguous B\"uchi Turing machine,\smallskip
 
\item the \emph{Gandy-Harrington topology} $\tau_{GH}$, generated by the set $\mathbb{B}_{GH}$ of languages accepted by some B\"uchi Turing machine.
\end{itemize}

 In \cite{SchwarzS10}, Schwarz and Staiger prove that $\tau_B$ and $\tau_A$ are metrizable. The topology $\tau_{GH}$ is second countable, $T_1$ and strong Choquet, but it is not regular and thus not metrizable and not Polish. However, there is a dense basic open set on which $\tau_{GH}$ is Polish and zero dimensional, which is sufficient in many applications. The topology $\tau_B$ is separable, by definition, because there are only countably many regular $\omega$-languages. It remains to see that it is completely metrizable to see that it is Polish. This is one of the main results proved in this paper.
 
\begin{them} \label{thm:main}
  Let $z \in \{ C,B,A,\delta\}$. Then $\tau_z$ is Polish and zero-dimensional.
\end{them}
 
 From this result, it is already possible to infer many properties of the space $\Sigma^{\mathbb{N}}$, equipped with the B\"uchi topology (see \cite{CFL17} for an extended list). In particular, we get some results about the $\sigma$-algebra generated by the $\omega$-regular languages. It is stratified in a hierarchy of length $\omega_1$ (the first uncountable
ordinal) and there are universal sets at each level of this
hierarchy. Notice that this $\sigma$-algebra coincides with the
$\sigma$-algebra of Borel sets for the Cantor topology. However, the levels of the Borel hierarchy differ for the two topologies. For instance, an $\omega$-regular set which is
non-${\bf \Pi}_2^0$ for the Cantor topology is clopen (i.e.,
${\bf \Delta}_1^0$) for the B\"uchi topology. Therefore the results about
the existence of universal sets at each level of the $\sigma$-algebra
generated by the $\omega$-regular languages are really new and
interesting. 

We also investigate, following a suggestion of H. Michalewski, whether it is possible to extend these results to the case of a space $T_\Sigma^\omega$ of infinite binary trees labelled with letters of the alphabet $\Sigma$. On the one hand, the automatic topology can be proved to be Polish in a similar way. On the other hand, we show that the B\"uchi topology (generated by the set of regular tree languages accepted by some B\"uchi tree automaton) and the Muller topology (generated by the set of regular tree languages accepted by some Muller tree automaton) are both non-Polish. However we prove that these two topologies have quite different properties: the first one is strong Choquet but not metrizable while the second one is metrizable but not strong Choquet. 

\section{Background}\indent

We first recall the notions required to understand fully the introduction
and the sequel (see for example
\cite{PerrinPin,Staiger97,Kechris94,Moschovakis09}).

\subsection{Theoretical computer science}\indent

A \emph{B\"uchi automaton} is a tuple
$\mathcal{A} = (\Sigma,Q,Q_i,Q_f,\delta )$, where $\Sigma$ is the input alphabet, $Q$ is the finite set of states, $Q_i$ and~$Q_f$ are the sets of initial and final states, and $\delta$ is the transition relation.  The transition relation~$\delta$ is a subset of $Q \times \Sigma \times Q$.

A \emph{run} on some sequence $\sigma \in \Sigma^{\mathbb{N}}$ is a
sequence $(q_n)_{i \in \mathbb{N}} \in Q^{\mathbb{N}}$ of states such that
$q_0$ is initial ($q_0 \in Q_i$) and $\big( q_i,\sigma(i),q_{i+1}\big)$ is a
transition in~$\delta$ for each $i \geq 0$.  It is \emph{accepting} if it visits
infinitely often final states, i.e., $q_i \in Q_f$ for infinitely many
$i$'s.  An input sequence~$\sigma$ is accepted if there exists an accepting
run on~$\alpha$.  The set of accepted inputs is denoted $L(\mathcal{A})$. A set of infinite words is called \emph{$\om$-regular} if it is equal to $L(\mathcal{A})$ for some automaton~$\mathcal{A}$ (see \cite{PerrinPin} for the basic notions about regular $\om$-languages, which are the $\om$-languages accepted by some B\"uchi or Muller automaton).

A B\"uchi automaton is actually similar to a \emph{classical finite automaton}.
A finite word $w$ of length~$n$ is accepted by some automaton~$\mathcal{A}$
if there is sequence $(q_i)_{i \leq n}$ of $n+1$ states such that $q_0$ is
initial ($q_0 \in Q_i$), $q_n$ is final ($q_n \in Q_f$) and
$(q_i,\sigma(i),q_{i+1})$ is a transition in~$\delta$ for each
$0 \leq i < n$.  The set of accepted finite words is denoted by
$U(\mathcal{A})$. A set of finite words is called \emph{regular} if it is
equal to $U(\mathcal{A})$ for some automaton~$\mathcal{A}$.

Let $U$ be a set of finite words, and $V$ be a set of finite or infinite words. We recall that 
$$U\cdot V=\{ u\cdot v\mid u\in U\wedge v\in V\} .$$ 
An infinite word is \emph{ultimately periodic} if it is of the form $u\cdot v^\omega$, where $u,v$ are finite words. The 
\emph{$\omega$-power} of a set~$U$ of finite words is defined by
\begin{displaymath}
  U^{\omega} = \{ \sigma \in \Sigma^{\mathbb{N}} \mid \exists \; 
               (w_i)_{i \in \mathbb{N}} \in U^{\mathbb{N}} \text{ such that } \sigma =
               w_0\cdot w_1\cdot w_2\cdots \}.
\end{displaymath}
The $\omega$-powers play a crucial role in the characterization of
$\omega$-regular languages (see \cite{Buchi62}).
 
\begin{thm}[B\"uchi] \label{buch} 
  Let $\Sigma$ be a finite alphabet, and $L\subseteq
  \Sigma^{\mathbb{N}}$. The following are equivalent:
\begin{enumerate}
\item $L$ is $\omega$-regular,
\item there are $2n$ regular languages $(U_i)_{i<n}$ and $(V_i)_{i<n}$ such
  that $L = \bigcup_{i<n}\ U_i\cdot V_i^\omega$.
\end{enumerate}
\end{thm}
 
In particular, each singleton $\{ uv^\omega\}$ formed by an ultimately periodic $\omega$-word is an $\omega$-regular language. We now recall some important properties of the class of $\omega$-regular languages (see \cite{PerrinPin} and~\cite{SchwarzS10}). Let $\Sigma$ be a set, and $\Sigma^*$ be the set of finite sequences of elements of $\Sigma$. If $w \in \Sigma^*$, then $w$ defines the usual basic clopen (i.e., closed and open) set
$N_w := \{\sigma \in \Sigma^{\mathbb{N}}\mid w\mbox{ is a prefix of
}\sigma\}$ of the Cantor topology $\tau_C$ (so
$\mathbb{B}_C := \{\emptyset\}\cup\{ N_w\mid w \in \Sigma^*\}$ is a basis
for $\tau_C$).
 
\begin{thm}[B\"uchi] \label{closure} 
  The class of $\omega$-regular languages contains the usual basic clopen
  sets and is closed under finite unions and intersections, taking
  complements, and projections (from a product alphabet onto one of its coordinates).
\end{thm}

 We now turn to the study of Turing machines (see \cite{CG78b, Staiger97}). A \emph{B\"uchi Turing machine} is a tuple $\mathcal{M} = (\Sigma ,\Gamma ,Q,q_0,Q_f,\delta )$, where $\Sigma$ and~$\Gamma$ are the input and tape alphabets satisfying $\Sigma \subseteq \Gamma$, $Q$ is the finite set of states, $q_0$ is the initial state, $Q_f$ is the set of
final states, and $\delta$ is the transition relation.  The relation~$\delta$ is a subset of $(Q \times \Gamma) \times (Q \times \Gamma \times \{ -1, 0, 1\} )$.

A \emph{configuration} of~$\mathcal{M}$ is a triple $(q,\gamma,j)$ where
$q \in Q$ is the current state, $\gamma \in \Gamma^{\mathbb{N}}$ is the
content of the tape and the non-negative integer $j \in \mathbb{N}$ is the
position of the head on the tape.

Two configurations $(q,\gamma,j)$ and $(q',\gamma',j')$ of~$\mathcal{M}$
are \emph{consecutive} if there exists a transition
$(q,a,q',b,d) \in \delta$ such that the following conditions are met.
\begin{enumerate}
\item $\gamma(j) = a$, $\gamma'(j) = b$ and $\gamma(i) = \gamma'(i)$ for
  each $i \neq j$.  This means that the symbol~$a$ is replaced by the 
  symbol~$b$ at the position~$j$ and that all the other symbols on the tape remain
  unchanged.
\item The two positions $j$ and~$j'$ satisfy the equality $j' = j+d$.
\end{enumerate}

A \emph{run} of the machine~$\mathcal{M}$ on some input
$\sigma \in \Sigma^{\mathbb{N}}$ is a sequence
$(p_i,\gamma_i,j_i)_{i\in\mathbb{N}}$ of consecutive configurations such
that $p_0 = q_0$, $\gamma_0 = \sigma$ and $j_0 = 0$.  The run is
\emph{accepting} if it visits infinitely often the final states, i.e.,
$p_i \in Q_f$ for infinitely many $i$'s.  The $\omega$-language accepted
by~$\mathcal{M}$ is the set of inputs~$\sigma$ such that there exists an
accepting run on~$\sigma$.

Notice that some other accepting conditions have been considered for the
acceptance of infinite words by Turing machines, like the 1' or Muller ones
(the latter one was firstly called 3-acceptance), see
\cite{CG78b,Staiger97}. Moreover, several types of required behaviour on the
input tape have been considered in the literature, see
\cite{Staiger99,FL-TM,Fin-ambTM}.

A \emph{B\"uchi automaton} $\mathcal{A}$ is in fact a B\"uchi Turing
machine whose head only moves forwards.  This means that each of its
transitions has the form $(p,a,q,b,d)$ where $d = 1$.  Note that the written
symbol~$b$ is never read.

\subsection{Descriptive set theory}\indent

 Classical descriptive set theory takes place in Polish topological spaces. We first recall that if $d$ is a distance on a set $X$, and $(x_n)_{n\in\mathbb{N}}$ is a sequence of elements of $X$, then the sequence  $(x_n)_{n\in\mathbb{N}}$ is called a \emph{Cauchy sequence} if  
$$\forall k \in \mathbb{N} ~~ \exists N \in \mathbb{N}~~  \forall p, p' \geq N ~~  d(x_p, x_{p'}) <  \frac{1}{2^k}.$$ 
\noindent In a topological space $X$ whose topology is induced by a distance $d$, the distance $d$ and the metric space $(X,d)$  are said to be {\it complete}  if every Cauchy sequence in $X$ is convergent. 
 
\begin{defi} A topological space $X$ is a \emph{Polish space} if it is
  \begin{enumerate}
  \item separable (there is a countable dense sequence $(x_n)_{n\in\mathbb{N}}$ in $X$),
  \item completely metrizable (there is a complete distance $d$ on $X$
    which is compatible with the topology of $X$).
  \end{enumerate}
\end{defi}

 The most classical hierarchy of topological complexity in descriptive set theory is the one given by the Borel classes. If $\bf\Gamma$ is a class of sets in metrizable spaces, then $\check {\bf\Gamma}\! :=\!\{\neg S\mid S\!\in\! {\bf\Gamma}\}$, and $(\bf\Gamma)_\sigma$ is the class of countable unions of sets in $\bf\Gamma$. Recall that the \emph{Borel hierarchy} is the inclusion from left to right in the following picture.\bigskip
 
\scalebox{0.74}{$$\!\!\!\!\!\!\!\!\!\!\!\!\!\!\!\!\!\!\!\!\!\!\!\xymatrix@1{ 
& & \boraone\! =\!\mbox{open} & & \boratwo\! =\! (\bormone )_\sigma & & & 
\boraxi\! =\! (\bigcup_{\eta <\xi}~\borme )_\sigma & \\ 
& \borone\! =\!\mbox{clopen} & & \bortwo\! =\!\boratwo\cap\bormtwo & & \cdots & \borxi\! =\!\boraxi\cap\bormxi & & \cdots\\ 
& & \bormone\! =\!\mbox{closed} & & \bormtwo\! =\!\check\boratwo & & & \bormxi\! =\!\check\boraxi & }$$}\bigskip

\noindent Above the Borel hierarchy sits the \emph{projective hierarchy}, which is the inclusion from left to right in the following picture.\bigskip
 
\scalebox{0.74}{$$\xymatrix@1{ 
& & \ana\! =\!\mbox{analytic} & & {\bf\Sigma}^1_2\! =\!\mbox{Projections of }\ca\mbox{ sets} & & & 
{\bf\Sigma}^1_{n+1}\! =\!\mbox{Projections of }{\bf\Pi}^1_n\mbox{ sets} & \\ 
& & \ca\! =\!\check\ana & & {\bf\Pi}^1_2\! =\!\check{\bf\Sigma}^1_2 & & & {\bf\Pi}^1_{n+1}\! =\!\check{\bf\Sigma}^1_{n+1} & }$$}\bigskip

\noindent Effective descriptive set theory is based on the notion of a recursive
function.  A function from $\mathbb{N}^k$ to $\mathbb{N}^l$ is said to be
\emph{recursive} if it is total and computable.  By extension, a relation
is called \emph{recursive} if its characteristic function is recursive.

\begin{defi}  
  A \emph{recursive presentation} of a Polish space $X$ is a pair
  $\big( (x_n)_{n\in\mathbb{N}},d\big)$ such that
  \begin{enumerate}
  \item $(x_n)_{n\in\mathbb{N}}$ is dense in $X$,
    
  \item $d$ is a compatible complete distance on $X$ such that the following
  relations $P$ and~$Q$ are recursive:
  \begin{align*}
    P(i,j,m,k) & \iff d(x_i,x_j) \leq \frac{m}{k+1}, \\
    Q(i,j,m,k) & \iff d(x_i,x_j) < \frac{m}{k+1}.
  \end{align*}
\end{enumerate}
A Polish space $X$ is \emph{recursively presented} if there is a recursive
presentation of it.
\end{defi}

 Note that the formula $(p,q) \mapsto 2^p(2q + 1) - 1$ defines a recursive
bijection $\mathbb{N}^2 \rightarrow \mathbb{N}$. One can check that the
coordinates of the inverse map are also recursive. They will be denoted
$n \mapsto (n)_0$ and $n \mapsto (n)_1$ in the sequel.  These maps will
help us to define some of the basic effective classes.

\begin{defi} \label{eff} 
  Let $\big( (x_n)_{n\in\mathbb{N}},d\big)$ be a recursive presentation of
  a Polish space $X$.

  \begin{enumerate}
  \item We fix a countable basis of $X$: $B(X,n)$ is the open ball
    $B_d(x_{(n)_0},\frac{( (n)_1)_0}{( (n)_1)_1+1})$.

  \item A subset $S$ of $X$ is \emph{semirecursive}, or \emph{effectively
      open} (denoted $S \in \Boraone$) if
    $$S = \bigcup_{n\in\mathbb{N}}{B\big( X,f(n)\big)}\mbox{,}$$ for some recursive
    function $f$.

  \item A subset $S$ of $X$ is \emph{effectively closed} (denoted
    $S \in \Bormone$) if its complement $\neg S$ is
    semirecursive.

  \item One can check that a product of two recursively presented Polish
    spaces has a recursive presentation, and that the Baire space
    $\mathbb{N}^\mathbb{N}$ has a recursive presentation.  A subset $S$ of
    $X$ is \emph{effectively analytic} (denoted $S \in \Ana$) if there
    is a $\Bormone$ subset $C$ of $X\times \mathbb{N}^\mathbb{N}$ such
    that
    \begin{displaymath}
      S =\pi_0[C] := \{ x \in X\mid\exists\alpha \in \mathbb{N}^\mathbb{N}
      \;\; (x,\alpha) \in  C\} .
    \end{displaymath}

\item A subset $S$ of $X$ is \emph{effectively co-analytic} (denoted $S \in \Ca$) if its
  complement $\neg S$ is effectively analytic, and
  \emph{effectively Borel} if it is in $\Ana$ and $\Ca$ (denoted
  $S \in \Borel$).

\item We will also use the following \emph{relativized classes}: if $X$,
  $Y$ are recursively presented Polish spaces and $y \in  Y$, then we say
  that $A\subseteq X$ is in $\Ana (y)$ if there is
  $S \in \Ana (Y\times X)$ such that $A = S_y := \{ x\in X\mid (y,x)\in S\}$. The class
  $\Ca (y)$ is defined similarly. We also set
  $\Borel (y) := \Ana (y)\cap\Ca (y)$.
\end{enumerate}
\end{defi}

The crucial link between the effective classes and the classical
corresponding classes is as follows: the class of analytic (resp.,
co-analytic, Borel) subsets of $Y$ is equal to
$\bigcup_{\alpha\in\mathbb{N}^\mathbb{N}}~\Ana (\alpha )$ (resp.,
$\bigcup_{\alpha\in\mathbb{N}^\mathbb{N}}~\Ca (\alpha )$,
$\bigcup_{\alpha\in\mathbb{N}^\mathbb{N}}~\Borel (\alpha )$). This allows to
use effective descriptive set theory to prove results of classical type. In
the sequel, when we consider an effective class in some
$\Sigma^{\mathbb{N}}$ with $\Sigma$ finite, we will always use a fixed
recursive presentation associated with the Cantor topology. The following
result is proved in \cite{Staiger99}, see also \cite{Fin-ambTM}.

\begin{thm} \label{TM} 
  Let $\Sigma$ be a finite alphabet, and $L\subseteq \Sigma^{\mathbb{N}}$.
  The following statements are equivalent:
  \begin{enumerate}
  \item $L = L(\mathcal{M})$ for some B\"uchi Turing machine $\mathcal{M}$,
  \item $L \in \Ana$.
  \end{enumerate}
\end{thm}

We now recall the \emph{strong Choquet game} played by two players on a
topological space~$X$.  Players 1 and~2 play alternatively.  At each
turn~$i$, Player 1 plays by choosing an open subset~$U_i$ and a
point~$x_i \in U_i$ such that $U_i \subseteq V_{i-1}$, where $V_{i-1}$ has
been chosen by Player 2 at the previous turn.  Player 2, plays by choosing
an open subset~$V_i$ such that $x_i \in V_i$ and $V_i \subseteq U_i$.
Player 2 wins the game if $\bigcap_{i\in\mathbb{N}}~V_i \neq \emptyset$.
We now recall some classical notions of topology.
 
\begin{defi} 
  A topological space~$X$ is said to be
  \begin{itemize}
  \item  $T_1$ if every singleton of $X$ is closed,
  \item \emph{regular} if for every point of $X$ and every open
    neighborhood $U$ of $x$, there is an open neighborhood $V$ of $x$ with
    $\overline{V}\subseteq U$,
  \item \emph{second countable} if its topology has a countable basis,
  \item \emph{zero-dimensional} if there is a basis made of clopen sets,
  \item \emph{strong Choquet} if $X$ is not empty and Player 2 has a
    winning strategy in the strong Choquet game. 
  \end{itemize}
\end{defi}

Note that every zero-dimensional space is regular. The following result is Theorem 8.18 in \cite{Kechris94}.
\begin{thm}[Choquet] \label{choquet} 
  A nonempty, second countable topological space is Polish if and only if
  it is $T_1$, regular, and strong Choquet.
\end{thm}

Let $X$ be a nonempty recursively presented Polish space. The
\emph{Gandy-Harrington topology} on $X$ is generated by the $\Ana$ subsets
of $X$, and denoted $\tau_{GH}^X$. By Theorem~\ref{TM}, this topology is
also related to automata and Turing machines. As there are some
effectively analytic sets whose complement is not analytic, the
Gandy-Harrington topology is not metrizable (in fact not regular) in
general (see 3E.9 in \cite{Moschovakis09}). In particular, it is not
Polish.

 Let $\bf\Gamma$ be a class of sets in Polish spaces. If $Y$ is a Polish space, then we say that $A\in {\bf\Gamma}(Y)$ is $\bf\Gamma$\emph{-complete} if, for each zero-dimensional Polish space $X$ and each $B\in {\bf\Gamma}(X)$, there is $f\! :\! X\!\rightarrow\! Y$ continuous such that $B\! =\! f^{-1}(A)$. By Section 22.B in \cite{Kechris94}, if 
$\bf\Gamma$ is of the form $\bf\Sigma$ or $\bf\Pi$ in the Borel or the projective hierarchy, and if $A$ is $\bf\Gamma$-complete, then $A$ is not in $\check {\bf\Gamma}$. Theorem 22.10 in \cite{Kechris94} gives a converse in the Borel hierarchy.

\section{Proof of Theorem~\ref{thm:main}}\indent

The proof of Theorem~\ref{thm:main} is organized as follows.  We provide
below four properties which ensure that a given topological space is strong
Choquet.  Then we use Theorem~\ref{choquet} to prove that the considered
spaces are indeed Polish.

Let $\Sigma$ be a countable alphabet. The set $\Sigma^{\mathbb{N}}$ is
equipped with the product topology of the discrete topology on~$\Sigma$,
unless another topology is specified.  This topology is induced by a
natural metric, called the \emph{prefix metric} which is defined as
follows. For $\sigma\neq\sigma' \in \Sigma^{\mathbb{N}}$, the distance~$d$
is given by
\begin{displaymath}
  d(\sigma,\sigma') = \frac{1}{2^r} 
  \quad\text{, where}\quad
  r = \min \{ n\in\mathbb{N} \mid \sigma(n) \neq \sigma'(n) \} .
\end{displaymath}
When $\Sigma$ is finite this topology is the classical Cantor topology.
When $\Sigma$ is countably infinite the topological space is homeomorphic
to the Baire space $\mathbb{N}^{\mathbb{N}}$.

Let $\Sigma$ and $\Gamma$ be two alphabets. The function which maps each
pair $(\sigma ,\gamma) \in \Sigma^{\mathbb{N}} \times \Gamma^{\mathbb{N}}$
to the element $\big(\sigma(0),\gamma(0)\big) ,\big(\sigma(1),\gamma(1)\big) ,\ldots$ of
$(\Sigma\times\Gamma)^{\mathbb{N}}$ is a homeomorphism between
$\Sigma^{\mathbb{N}}\times \Gamma^{\mathbb{N}}$ and
$(\Sigma\times\Gamma)^{\mathbb{N}}$ allowing us to identify these two
spaces.

If $\Sigma$ is a set, $\sigma \in \Sigma^{\mathbb{N}}$ and
$l \in \mathbb{N}$, then $\sigma\vert l$ is the prefix of $\sigma$ of
length~$l$.

We set $2 := \{ 0,1\}$ and
$\mathbb{P}_\infty := \{\alpha \in 2^{\mathbb{N}}\mid\forall k \in
\mathbb{N} \;\; \exists i \geq k \;\; \alpha (i) = 1\}$.
This latter set is simply the set of infinite words over the alphabet $2$ having infinitely many $1$'s.

We will work in the spaces of the form $\Sigma^{\mathbb{N}}$, where
$\Sigma$ is a finite set with at least two elements. We consider a topology
$\tau_\Sigma$ on $\Sigma^{\mathbb{N}}$, and a basis $\mathbb{B}_\Sigma$ for
$\tau_\Sigma$.  We consider the following properties of the family
$(\tau_\Sigma ,\mathbb{B}_\Sigma )_\Sigma$, using the previous
identification of $\Sigma^{\mathbb{N}}\times \Gamma^{\mathbb{N}}$ and
$(\Sigma\times\Gamma)^{\mathbb{N}}$:

\begin{enumerate}[label=(P\arabic*)]
\item $\mathbb{B}_\Sigma$ contains the usual basic clopen sets
  $N_w$,\smallskip

\item $\mathbb{B}_\Sigma$ is closed under finite unions and
  intersections,\smallskip

\item $\mathbb{B}_\Sigma$ is closed under projections, in the sense that
  if $\Gamma$ is a finite set with at least two elements and
  $L \in \mathbb{B}_{\Sigma\times\Gamma}$, then
  $\pi_0[L] \in \mathbb{B}_\Sigma$,\smallskip

\item for each $L \in \mathbb{B}_\Sigma$ there is a closed subset
  $C$ of $\Sigma^{\mathbb{N}} \times \mathbb{P}_\infty$ (i.e., $C$ is the
  intersection of a $\tau_C$-closed subset of the Cantor space
  $\Sigma^{\mathbb{N}} \times 2^{\mathbb{N}}$ with
  $\Sigma^{\mathbb{N}} \times \mathbb{P}_\infty$), which is in
  $\mathbb{B}_{\Sigma\times 2}$, and such that $L =\pi_0[C]$.
\end{enumerate}

\begin{thm} \label{main}  
  Assume that a family $(\tau_\Sigma ,\mathbb{B}_\Sigma )_\Sigma$
  satisfies Properties (P1)-(P4). Then the topologies $\tau_\Sigma$ are
  strong Choquet.
\end{thm}

\noindent\bf Proof.\rm\ 
  We first describe a strategy $\tau$ for Player~2. Player~1 first plays
  $\sigma_0 \in \Sigma^{\mathbb{N}}$ and a $\tau_\Sigma$-open neighborhood
  $U_0$ of $\sigma_0$. Let $L_0$ in $\mathbb{B}_\Sigma$ with
  $\sigma_0 \in L_0\subseteq U_0$. Property (P4) gives $C_0$ with
  $L_0 =\pi_0[C_0]$. This gives $\alpha_0 \in \mathbb{P}_\infty$ such that
  $(\sigma_0,\alpha_0) \in C_0$. We choose $l^0_0 \in \mathbb{N}$ big
  enough to ensure that if 
  $$s^0_0 := \alpha_0\vert l^0_0\mbox{,}$$ 
  then $s^0_0$ has at least a coordinate equal to $1$. We set $w_0 := \sigma_0\vert 1$ and
  $V_0 := \pi_0[C_0\cap (N_{w_0}\times N_{s^0_0})]$. By Properties
  (P1)-(P3), $V_0$ is in $\mathbb{B}_\Sigma$ and thus
  $\tau_\Sigma$-open. Moreover,
  $\sigma_0 \in V_0\subseteq L_0\subseteq U_0$, so that Player~2 respects
  the rules of the game if he plays $V_0$.

  Now Player~1 plays $\sigma_1 \in V_0$ and a $\tau_\Sigma$-open
  neighborhood $U_1$ of $\sigma_1$ contained in $V_0$. Let $L_1$ in
  $\mathbb{B}_\Sigma$ with $\sigma_1 \in L_1\subseteq U_1$. Property (P4)
  gives $C_1$ with $L_1 =\pi_0[C_1]$. This gives
  $\alpha_1 \in \mathbb{P}_\infty$ such that $(\sigma_1,\alpha_1) \in
  C_1$. We choose $l^1_0 \in \mathbb{N}$ big enough to ensure that if
  $s^1_0 := \alpha_1\vert l^1_0$, then $s^1_0$ has at least one coordinate
  equal to $1$. As $\sigma_1 \in V_0$, there is
  $\alpha'_0 \in \mathbb{P}_\infty$ such that
  $(\sigma_1,\alpha'_0) \in C_0\cap (N_{w_0}\times N_{s^0_0})$. We choose
  $l^0_1 > l^0_0$ big enough to ensure that if
  $s^0_1 := \alpha'_0\vert l^0_1$, then $s^0_1$ has at least two
  coordinates equal to $1$. We set $w_1 := \sigma_1\vert 2$ and
  $V_1 := \pi_0[C_0\cap (N_{w_1}\times N_{s^0_1})]\cap\pi_0[C_1\cap
  (N_{w_0}\times N_{s^1_0})]$. Here again, $V_1$ is
  $\tau_\Sigma$-open. Moreover, $\sigma_1 \in V_1\subseteq U_1$ and Player
  2 can play $V_1$.

  Next, Player~1 plays $\sigma_2 \in V_1$ and a $\tau_\Sigma$-open
  neighborhood $U_2$ of $\sigma_2$ contained in $V_1$. Let $L_2$ in
  $\mathbb{B}_\Sigma$ with $\sigma_2 \in L_2\subseteq U_2$. Property (P4)
  gives $C_2$ with $L_2 =\pi_0[C_2]$. This gives
  $\alpha_2 \in \mathbb{P}_\infty$ such that $(\sigma_2,\alpha_2) \in
  C_2$. We choose $l^2_0 \in \mathbb{N}$ big enough to ensure that if
  $s^2_0 := \alpha_2\vert l^2_0$, then $s^2_0$ has at least one coordinate
  equal to $1$. As $\sigma_2 \in V_1$, there is
  $\alpha'_1 \in \mathbb{P}_\infty$ such that
  $(\sigma_2,\alpha'_1) \in C_1\cap (N_{w_0}\times N_{s^1_0})$. We choose
  $l^1_1 > l^1_0$ big enough to ensure that if
  $s^1_1 := \alpha'_1\vert l^1_1$, then $s^1_1$ has at least two
  coordinates equal to $1$. As $\sigma_2 \in V_1$, there is
  $\alpha''_0 \in \mathbb{P}_\infty$ such that
  $(\sigma_2,\alpha''_0) \in C_0\cap (N_{w_1}\times N_{s^0_1})$. We choose
  $l^0_2 > l^0_1$ big enough to ensure that if
  $s^0_2 := \alpha''_0\vert l^0_2$, then $s^0_2$ has at least three
  coordinates equal to $1$. We set $w_2 := \sigma_2\vert 3$ and
  $V_2 := \pi_0[C_0\cap (N_{w_2}\times N_{s^0_2})]\cap\pi_0[C_1\cap
  (N_{w_1}\times N_{s^1_1})]\cap\pi_0[C_2\cap (N_{w_0}\times
  N_{s^2_0})]$. Here again, $V_2$ is $\tau_\Sigma$-open. Moreover,
  $\sigma_2 \in V_2\subseteq U_2$ and Player~2 can play $V_2$.

  If we go on like this, we build $w_l \in \Sigma^{l+1}$ and
  $s^n_l \in 2^*$ such that $w_0\subseteq w_1\subseteq ...$ and
  $s^n_0\!\subsetneqq\! s^n_1\!\subsetneqq\! ...$ This allows us to define
  $\sigma := \mbox{lim}_{l\rightarrow\infty}~w_l \in \Sigma^{\mathbb{N}}$
  and, for each $n \in \mathbb{N}$,
  $\alpha_n :=\mbox{lim}_{l\rightarrow\infty}~s^n_l \in
  2^{\mathbb{N}}$. Note that $\alpha_n \in \mathbb{P}_\infty$ since $s^n_l$
  has at least $l + 1$ coordinates equal to $1$. As
  $(\sigma ,\alpha_n)$ is the limit of $(w_l,s^n_l)$ as $l$ goes to
  infinity and $N_{w_l}\times N_{s^n_l}$ meets $C_n$ (which is closed in
  $\Sigma^{\mathbb{N}}\times \mathbb{P}_\infty$ in the sense of Property (P4)),
  $(\sigma ,\alpha_n) \in C_n$. Thus
$$\sigma \in \bigcap_{n\in\mathbb{N}}~\pi_0[C_n] =
\bigcap_{n\in\mathbb{N}}~L_n\subseteq \bigcap_{n\in\mathbb{N}}~U_n\subseteq
\bigcap_{n\in\mathbb{N}}~V_n\mbox{,}$$ so that $\tau$ is winning for Player
2.\hfill{$\square$}

\subsection{The Gandy-Harrington topology}\indent

We have already mentioned the fact that the Gandy-Harrington topology is not Polish
in general. However, it is almost Polish since it fulfills Properties
(P1)-(P4).

Let $\Sigma$ be a finite alphabet with at least two elements and $X$ be
the space $\Sigma^{\mathbb{N}}$ equipped with the topology
$\tau_\Sigma :=\tau_{GH}^X$ generated by the family~$\mathbb{B}_\Sigma$
of $\Ana$ subsets of $X$. Note that the assumption of Theorem~\ref{main}
are satisfied. Indeed, (P1)-(P3) come from 3E.2 in
\cite{Moschovakis09}. For (P4), let $F$ be a $\Bormone$ subset of
$X\times \mathbb{N}^\mathbb{N}$ such that $L =\pi_0[F]$.  Let $\varphi$ be
the function from $\mathbb{N}^\mathbb{N}$ to $2^\mathbb{N}$ defined by
$$\varphi(\beta) = 0^{\beta(0)}10^{\beta (1)}1\ldots$$  
Note that $\varphi$ is a homeomorphism from $\mathbb{N}^\mathbb{N}$ onto $\mathbb{P}_\infty$,
and recursive (which means that the relation
$\varphi(\beta) \in N(2^{\mathbb{N}},n)$ is semirecursive in $\beta$ and
$n$). This implies that $C := (\operatorname{Id} \times \varphi )[F]$ is
suitable (see 3E.2 in \cite{Moschovakis09}).

Note that $\tau_\Sigma$ is second countable since there are only countably
many $\Ana$ subsets of $X$ (see 3F.6 in \cite{Moschovakis09}), $T_1$ since
it is finer than the usual topology by the property (P1), and strong
Choquet by Theorem~\ref{main}.

One can show that there is a dense basic open subset $\Omega_X$ of
$(X,\tau_\Sigma)$ such that $S \cap \Omega_X$ is a clopen subset of
$(\Omega_X ,\tau_\Sigma )$ for each $\Ana$ subset~$S$ of $X$ (see
\cite{Lec09a}). In particular, $(\Omega_X,\tau_\Sigma )$ is
zero-dimensional, and regular. As it is, just like $(X,\tau_\Sigma)$,
second countable, $T_1$ and strong Choquet, $(\Omega_X,\tau_\Sigma )$
is a Polish space, by Theorem~\ref{choquet}.

\subsection{The B\"uchi topology}\indent
 
Let $\Sigma$ be a finite alphabet with at least two symbols, and $X$ be
the space $\Sigma^{\mathbb{N}}$ equipped with the B\"uchi topology $\tau_B$
generated by the family $\mathbb{B}_B$ of $\omega$-regular languages in
$X$. Theorem 29 in \cite{SchwarzS10} shows that $\tau_B$ is metrizable. We now give a distance which is compatible with $\tau_B$. This metric was used in \cite{HoffmannSS16} (Theorem 2 and Lemma 21 and several corollaries following Lemma 21). A similar argument for subword metrics is in Section 4 in \cite{HoffmannS15}. If 
$\mathcal{A}$ is a B\"uchi automaton, then we denote $|\mathcal{A}|$ the number of states of 
$\mathcal{A}$. We say that a B\"uchi automaton {\it separates} $x$ and $y$ if and only if 
$$\big( x \in  L(\mathcal{A})\wedge y\!\notin\! L(\mathcal{A})\big)\vee
\big( y \in  L(\mathcal{A})\wedge x\!\notin\! L(\mathcal{A})\big) .$$ 
The distance $\delta$ on $\Sigma^{\mathbb{N}}$ is then defined as follows:
\[
\delta(x,y) =\left\{\!\!\!\!\!\!\!
\begin{array}{ll}
& 0\mbox{ if }x = y\mbox{,}\cr
& \frac{1}{2^n}\mbox{ if }x \neq y\mbox{,}
\end{array}
\right.
\quad\text{for }x,y \in \Sigma^{\mathbb{N}} .
\]
where $n := \mbox{min}\{ |\mathcal{A}|\mid\mathcal{A}\text{ is a B\"uchi automaton which separates }x\text{ and }y\}$. We now describe some properties of the map $\delta$. This is the occasion to illustrate the notion of a complete metric.

\begin{prop} \label{propo} The following properties of $\delta$ hold:
\begin{enumerate}
  \item the map $\delta$ defines a distance on $\Sigma^{\mathbb{N}}$,
  \item the distance $\delta$ is compatible with $\tau_B$,
  \item the distance $\delta$ is not complete.
  \end{enumerate}
\end{prop}

\noindent\bf Proof.\rm\  1. If $x,y \in \Sigma^{\mathbb{N}}$,
then $\delta(x,y) =\delta(y,x)$, by the definition of $\delta$. Let
$x,y,z \in \Sigma^{\mathbb{N}}$, and assume that
$\delta(x,y) +\delta(y,z) <\delta(x,z) = \frac{1}{2^n}$. Then
$\delta(x,y) < \frac{1}{2^n}$ and $\delta(y,z) < \frac{1}{2^n}$ hold. In
particular, if $\mathcal{A}$ is a B\"uchi automaton with $n$ states then it
does not separate $x$ and $y$ and similarly it does not separate $y$ and
$z$. Thus either $x,y,z \in L(\mathcal{A})$ or
$x, y, z \notin L(\mathcal{A})$. This implies that the B\"uchi automaton
$\mathcal{A}$ does not separate $x$ and $z$. As this holds for every
B\"uchi automaton with $n$ states, $\delta(x,z) < \frac{1}{2^n}$. This
leads to a contradiction and thus
$\delta(x,z) \leq \delta(x,y) +\delta(y,z)$ for all
$x,y,z \in \Sigma^{\mathbb{N}}$. This shows that $\delta$ is a distance on
$\Sigma^{\mathbb{N}}$.

2. Recall that an open set for this topology is a union of $\om$-languages accepted by some 
B\"uchi automaton. Let then $L(\mathcal{A})$ be an $\om$-language accepted
by some B\"uchi automaton $\mathcal{A}$ having $n$ states, and
$x\in L(\mathcal{A})$. We now show that the open ball $B(x, \frac{1}{2^{n+1}})$
with center $x$ and $\delta$-radius $\frac{1}{2^{n+1}}$ is a subset of
$L(\mathcal{A})$. Indeed, if $\delta(x,y) < \frac{1}{2^{n+1}} < \frac{1}{2^n}$,
then $x$ and $y$ cannot be separated by any B\"uchi automaton with $n$
states, and thus $y\in L(\mathcal{A})$. This shows that $L(\mathcal{A})$
(and therefore any open set for $\tau_B$) is open for the topology induced
by the distance $\delta$. Conversely, let $B(x, r)$ be an open ball for the
distance $\delta$, where $r > 0$ is a positive real. It is clear from
the definition of the distance $\delta$ that we may only consider the case
$r = \frac{1}{2^n}$ for some natural number $n$. Then $y \in B(x, \frac{1}{2^n})$ if and
only if $x$ and $y$ cannot be separated by any B\"uchi automaton with
$p \leq  n$ states. Therefore the open ball $B(x, \frac{1}{2^n})$ is the
intersection of the regular $\om$-languages $L(\mathcal{A}_i)$ for some B\"uchi
automata $\mathcal{A}_i$ having $p \leq  n$ states and such that
$x \in L(\mathcal{A}_i)$, and of the regular $\om$-languages
$\Sigma^{\mathbb{N}}\!\setminus\! L(\mathcal{B}_i)$ for some B\"uchi automata
$\mathcal{B}_i$ having $p \leq  n$ states and such that
$x\!\notin\! L(\mathcal{B}_i)$. The class of regular $\om$-languages being
closed under taking complements and finite intersections, the open ball
$B(x, \frac{1}{2^n})$ is actually a regular $\om$-language and thus an open set
for $\tau_B$.

3. Without loss of generality, we
set $\Sigma =2$ and we consider, for a natural number $n \geq  1$,
the $\om$-word $X_n = 0^{n!}\cdot 1 \cdot 0^\om$ over the alphabet $2$ having only one symbol $1$ after $n!$ symbols~$0$, where
$n! := n\times (n - 1)\times \cdots\times 2\times 1$. Let now
$m > n > k$ and $\mathcal{A}$ be a B\"uchi automaton with $k$
states. Using a classical pumping argument, we can see that the automaton
$\mathcal{A}$ cannot separate $X_n$ and $X_m$. Indeed, assume first that
$X_n \in L(\mathcal{A})$. Then, when reading the first $k$ symbols $0$ of
$X_n$, the automaton enters at least twice in a same state $q$. This
implies that: $(\exists p \leq k) ~~~   (\forall l\geq 1) ~~~  0^{n!+lp}\cdot 1 \cdot 0^\om \in L(\mathcal{A})$.
\nl In particular
$m! = n!\times (n + 1)\times \cdots\times m = n! + n!\times
\Big(\big( (n + 1)\times \cdots\!\times m\big) - 1\Big)$ is of this
form and thus $X_m \in L(\mathcal{A})$. A very similar pumping argument
shows that if $X_m \in L(\mathcal{A})$, then $X_n \in L(\mathcal{A})$. This
shows that $\delta(X_n,X_m) <  \frac{1}{2^k}$ and finally that the sequence
$(X_n)$ is a Cauchy sequence for the distance $\delta$. On the other hand
if this sequence was converging to an $\om$-word $x$ then $x$ should be the
word $0^\omega$ because $\tau_B$ is finer than $\tau_C$. But $0^\omega$ is
an ultimately periodic word and thus it is an isolated point for
$\tau_B$. This leads to a contradiction, and thus the distance
$\delta$ is not complete because the sequence $(X_n)$ is a Cauchy sequence
which is not convergent.\hfill{$\square$}\bigskip

 Proposition~\ref{propo} gives a motivation for deriving Theorem~\ref{thm:main} from Theorem~\ref{main}. Note that the assumption of Theorem~\ref{main} are satisfied. Indeed, (P1)-(P3) come from Theorem~\ref{closure}. We now check (P4).
 
\begin{lem} \label{lemm} 
  Let $\Sigma$ be a finite set with at least two elements, and
  $L\subseteq \Sigma^{\mathbb{N}}$ be an $\omega$-regular language. Then
  there is a closed subset $C$ of
  $\Sigma^{\mathbb{N}}\times\mathbb{P}_\infty$, which is $\omega$-regular
  as a subset of $(\Sigma\times 2)^{\mathbb{N}}$ identified with
  $\Sigma^{\mathbb{N}}\times 2^{\mathbb{N}}$, and such that $L =\pi_0[C]$.
\end{lem}

\noindent\bf Proof.\rm\ 
  Let $\mathcal{A} = (\Sigma, Q, \delta, Q_i, Q_f)$ be a B\"uchi automaton
  and let $L = L(\mathcal{A})$ be its set of accepted words.  Let $\chi_f$
  be the characteristic function of $Q_f$.  It maps the 
  state~$q$ to $1$ if $q \in Q_f$, and to~$0$ otherwise.  The
  function~$\chi_f$ is extended to $Q^{\mathbb{N}}$ by setting
  $\alpha = \chi_f((q_n)_{n\in \mathbb{N}})$ where $\alpha(n) = \chi_f(q_n)$.
  Note that a run~$\rho$ of~$\mathcal{A}$ is accepting if and only if
  $\chi_f(\rho) \in \mathbb{P}_\infty$. Let $C$ be the subset of $\Sigma^{\mathbb{N}}\times
  \mathbb{P}_\infty$ defined by
  \begin{displaymath}
    C := \big\{ (\sigma,\alpha) \in \Sigma^{\mathbb{N}} \times \mathbb{P}_\infty
         \mid\exists \rho 
         \text{ run of $\mathcal{A}$ on $\sigma$ \text{such that} 
                $\alpha = \chi_f(\rho)$}
         \big\}.
  \end{displaymath}
By the definition of $C$, $L =\pi_0[C]$.  Let $K$ be the subset of
  $\Sigma^{\mathbb{N}}\times 2^{\mathbb{N}} \times Q^{\mathbb{N}}$ 
  defined by
  \begin{displaymath}
    K := \big\{ (\sigma,\alpha,\rho) \in 
                \Sigma^{\mathbb{N}} \times 2^{\mathbb{N}}\times Q^{\mathbb{N}}
         \mid \rho 
         \text{ is a run of $\mathcal{A}$ on $\sigma$ \text{such that}
                $\alpha = \chi_f(\rho)$}
         \big\}.
  \end{displaymath}
  As $K$ is compact as a closed subset of a compact space and
  $C =\pi_{\Sigma^{\mathbb{N}}\times 2^{\mathbb{N}}}[K]\cap (\Sigma^{\mathbb{N}}\times \mathbb{P}_\infty )$, $C$ is a closed subset of
  $\Sigma^{\mathbb{N}}\times \mathbb{P}_\infty$. It remains to show that
  $C$ is indeed $\omega$-regular.  Let $\Delta$ be defined by
  \begin{displaymath}
    \Delta := \big\{\big( p,(a,\varepsilon ),q\big) \in  
              Q\times (\Sigma\times 2)\times Q\mid 
              (p,a,q) \in \delta \wedge 
              (\varepsilon = 1\iff p \in Q_f)\big\} .
  \end{displaymath}
This allows us to define a B\"uchi automaton by
 $\mathcal{A}' := (\Sigma\times 2,Q,\Delta ,Q_i,Q_f)$. Note that
$$\begin{array}{ll}
    (\sigma ,\alpha ) \in  L(\mathcal{A}')\!\!\!\!
    & \Leftrightarrow \exists (s_i)_{i\in\mathbb{N}} \in  Q^{\mathbb{N}}~~
      \big( s_0 \in  Q_i\ \wedge\ \forall i \in \mathbb{N}~~(s_i,\big(\sigma (i),\alpha (i)\big),s_{i+1}) \in \Delta\big)\ \wedge\cr 
    & \hfill{\forall k \in \mathbb{N}~~\exists i \geq  k~~s_i \in  Q_f)}\cr
    & \Leftrightarrow\alpha \in \mathbb{P}_\infty\ \wedge\ \exists (s_i)_{i\in\mathbb{N}} \in  Q^{\mathbb{N}}~~
      \big( s_0 \in  Q_i\ \wedge\ \forall i \in \mathbb{N}~~(s_i,\sigma (i),s_{i+1}) \in \delta\ \wedge\cr 
    & \hfill{(\alpha (i) = 1\Leftrightarrow s_i \in  Q_f)\big)}\cr
    & \Leftrightarrow (\sigma ,\alpha ) \in  C.\cr
  \end{array}$$
  Thus $C = L(\mathcal{A}')$ is $\omega$-regular.
\hfill{$\square$}

\begin{cor} \label{buc} Let $\Sigma$ be a finite set with at least two elements. Then the B\"uchi topology $\tau_B$ is zero dimensional and Polish.\end{cor}

\noindent\bf Proof.\rm\ 
  As there are only countably many possible automata (up to
  identifications), $\mathbb{B}_B$ is countable. This shows that $\tau_B$
  is second countable. It is $T_1$ since it is finer than the usual
  topology by Property (P1), and strong Choquet by Theorem
  \ref{main}. Moreover, it is zero-dimensional since the class of
  $\omega$-regular languages is closed under taking complements (see
  Theorem \ref{closure}). It remains to apply Theorem~\ref{choquet}.
\hfill{$\square$}

\subsection{The other topologies}

\begin{lem} \label{lemfer} Let $(X,\tau )$ be a Polish space, and $(C_n)_{n\in\mathbb{N}}$ be a sequence of closed subsets of $(X,\tau )$. Then the topology generated by $\tau\cup\{ C_n\mid n \in \mathbb{N}\}$ is Polish.\end{lem}

\noindent\bf Proof.\rm\ By Lemma 13.2 in \cite{Kechris94}, the topology $\tau_n$ generated by
  $\tau\cup\{ C_n\}$ is Polish. By Lemma 13.3 in \cite{Kechris94}, the topology
  $\tau_\infty$ generated by $\bigcup_{n\in\mathbb{N}}~\tau_n$ is
  Polish. Thus the topology generated by
  $\tau\cup\{ C_n\mid n \in \mathbb{N}\}$, which is $\tau_\infty$, is
  Polish.\hfill{$\square$}\bigskip

\noindent\bf Proof of Theorem~\ref{thm:main}.\rm\ It is well known that $(\Sigma^{\mathbb{N}},\tau_C)$ is metrizable and compact, and thus Polish, and zero-dimensional.
\begin{itemize}

\item By Theorem 3.4 in \cite{Louveau80}, the implication (iii)
  $\Rightarrow$ (i), $\Borel (\Sigma^{\mathbb{N}})$ is a basis for a
  zero-dimensional Polish topology on $\Sigma^{\mathbb{N}}$. Recall that a
  B\"uchi Turing machine is {\it unambiguous} if every $\omega$-word
  $\sigma \in \Sigma^{\mathbb{N}}$ has at most one accepting run. By
  Theorem 3.6 in \cite{Fin-ambTM}, a subset of $\Sigma^{\mathbb{N}}$ is
  $\Borel$ if and only if it is accepted by some unambiguous B\"uchi Turing
  machine. Therefore $\mathbb{B}_\delta =\Borel (\Sigma^{\mathbb{N}})$ is a
  basis for the zero-dimensional Polish topology $\tau_\delta$.

\item Corollary \ref{buc} gives the result for the B\"uchi topology.
 
\item Lemma \ref{lemfer} shows that the automatic topology is Polish since it refines the usual product topology on $\Sigma^{\mathbb{N}}$. For this reason also, it is zero-dimensional.
\hfill{$\square$}
\end{itemize}

\section{The B\"uchi and Muller topologies on a  space of trees}\indent

The notion of a B\"uchi automaton has been extended to the case of a  B\"uchi tree automaton reading infinite binary trees whose nodes are labelled by letters of  a finite alphabet. 
We now recall this notion and some related ones. 

 A node of an infinite binary tree is represented by a finite  word over 
the alphabet $\{l, r\}$ where $l$ means ``left" and $r$ means ``right". An 
infinite binary tree whose nodes are labelled  in $\Si$ is identified with a function
$t: \{l, r\}^\star \ra \Si$. The set of  infinite binary trees labelled in $\Si$ will be 
denoted $T_\Si^\om$.

 A finite binary tree is like an ``initial finite subtree'' of an  infinite binary tree. Thus it can be represented by a function 
$s: S \subseteq  \{l, r\}^\star \ra \Si$, where $S$ is a finite subset of $ \{l, r\}^\star$ which is closed under prefix. 
If $t \in T_\Si^\om$ is an infinite binary tree, and $n\geq 0$ is an integer, then we denote by $t | n$ the initial finite subtree of $t$ whose domain is equal to $\{l, r\}^{\leq n}$, where 
 $\{l, r\}^{\leq n}$ is the set of finite words over the alphabet $\{l, r\}$ of length smaller than or equal to $n$. 

 Let $t$ be an infinite binary  tree. A branch $B$ of $t$ is a subset of the set of nodes of $t$ which 
is linearly ordered by the tree partial order $\sqsubseteq$ and which 
is closed under prefix relation (i.e., if  $x$ and $y$ are nodes of $t$ such that $y\in B$ and $x \sqsubseteq y$, then $x\in B$). A branch $B$ of a tree is said to be maximal if and only if there is no other branch of $t$ which strictly contains $B$. Let $t$ be an infinite binary tree in $T_\Si^\om$. If $B$ is a maximal branch of $t$,
then this branch is infinite. Let $(u_i)_{i\geq 0}$ be the enumeration of the nodes in $B$
which is strictly increasing for the prefix order. 
  The infinite sequence of the labels of the nodes of  such a maximal 
branch $B$, i.e., $t(u_0)t(u_1) \cdots t(u_n) \cdots $, is called a path. It is an $\om$-word 
over the alphabet $\Si$.

  Let then $L\subseteq \Si^\om$ be an $\omega$-language over $\Si$. We denote $\exists \mathrm{Path}(L)$  the set of 
infinite trees $t$ in $T_\Si^\om$ such that $t$ has (at least) one path in $L$.

 We now define the tree automata and the recognizable tree languages. 

\begin{defi} A (non deterministic) tree automaton  is a quadruple $\mathcal{A}=(\Si,Q,q_0,\Delta )$, where $\Sigma$ is the finite input alphabet, $Q$ is the finite set of states, $q_0 \in Q$ is the initial state and $\Delta \subseteq  Q \times   \Si   \times  Q \times   Q$ is the transition relation. The tree automaton  $\mathcal{A}$ is said to be 
deterministic if the relation  $\Delta$ is a functional one, i.e.,  if  for each $(q, a) \in Q \times   \Si $ there is at most one pair of states $(q', q'')$ 
such that  $(q, a, q', q'') \in \Delta$. 

 A run of the tree automaton  $\mathcal{A}$ on an infinite binary tree $t\in T_\Si^\om$ is an infinite binary tree $\rho \in T_Q^\om$ such that

 \begin{enumerate}[label=(\alph*)]
 \item  $\rho (\lambda)=q_0$,
 
 \item for each $u \in \{l, r\}^\star$,  $\big(\rho(u), t(u), \rho(u.l), \rho(u.r)\big)  \in \Delta$. 
\end{enumerate}
\end{defi}

\begin{defi}
A B\"uchi (non deterministic) tree automaton  is a  tuple 
$$\mathcal{A}=(\Si,Q,q_0,Q_f,\Delta )\mbox{,}$$ 
where $(\Si,Q,q_0,\Delta )$ is a tree automaton and $Q_f\subseteq Q$ is the set of accepting states.

 A run $\rho$ of the  B\"uchi tree automaton $\mathcal{A}$ on an infinite binary tree $t\in T_\Si^\om$ is said to be accepting if for each path of $\rho$ there is some accepting state appearing infinitely  often on this path. 

 The tree language $L(\mathcal{A})$ accepted by the B\"uchi tree automaton $\mathcal{A}$ is the set of infinite binary trees $t\in T_\Si^\om$ 
such that there is (at least) one accepting run of $\mathcal{A}$ on $t$. 
\end{defi}

\begin{defi}
A Muller  (non deterministic) tree automaton  is a  tuple 
$$\mathcal{A}=(\Si,Q,q_0,Q_f,\Delta )\mbox{,}$$ 
where $(\Si,Q,q_0,\Delta )$ is a tree automaton and $\mathcal{F} \subseteq 2^Q$ is the collection of  designated state sets. 

 A run $\rho$ of the  Muller  tree automaton $\mathcal{A}$ on an infinite binary tree $t\in T_\Si^\om$ is said to be accepting if  
for each path $p$ of $\rho$, the set of   states appearing infinitely  often on this path is in $\mathcal{F}$. 

 The tree language $L(\mathcal{A})$ accepted by the  Muller tree automaton $\mathcal{A}$ is the set of infinite binary trees $t\in T_\Si^\om$ 
such that there is (at least) one accepting run of $\mathcal{A}$ on $t$. 

 The class $REG$ of regular, or recognizable, tree languages is the class of  tree  languages accepted by some Muller automaton. 
\end{defi}

\begin{Rem}
Each tree language accepted by some (deterministic)  B\"uchi automaton is also accepted by some (deterministic) 
Muller automaton. A tree language is accepted by some Muller tree  automaton if and only if it is accepted by some Rabin tree automaton. We refer for instance to 
\cite{Thomas90,PerrinPin} for the definition of a Rabin tree  automaton. 
\end{Rem}

\begin{Exa}
Let  $L\subseteq \Si^\om$ be a regular $\om$-language. Then the set  $\exists \mathrm{Path}(L) \subseteq  T_\Si^\om$ is accepted by some  
B\"uchi tree automaton, hence also by some Muller  tree automaton. 

 The set of infinite binary trees $t\in T_\Si^\om$ having all their paths in $L$, denoted $\fa  \mathrm{Path}(L)$, 
  is accepted by some deterministic Muller  tree automaton. It is in fact the 
complement of the set $\exists \mathrm{Path}(\Sio-L) $. 
\end{Exa}

 There is a natural topology on the set $T_\Si^\om$ \cite{Moschovakis80,LescowThomas,Kechris94}. 
It is defined by the following distance. Let $t$ and $s$ be two distinct infinite trees in $T_\Si^\om$. 
Then the distance between $t$ and $s$ is $\frac{1}{2^n}$, where $n$ is the smallest integer 
such that $t(x)\neq s(x)$ for some word $x\in \{l, r\}^\star$ of length $n$.

 Let $T_0$ be a set of finite labelled trees, and $T_0\cdot T_\Si^\om$ be the set of infinite binary trees which extend some finite labelled binary tree $t_0\in T_0$. Here, $t_0$ is here a sort of prefix, an ``initial subtree" of a tree in $t_0\cdot T_\Si^\om$. The open sets are then of the form $T_0\cdot T_\Si^\om$. 

 It is well known that the set $T_\Si^\om$, equipped with this topology, is homeomorphic to the Cantor set $2^\om$, hence  also to the topological spaces 
$\Sio$, where $\Si$ is a finite alphabet having at least two letters. 

 We are going to use some notation similar to the one used in the case of the space $\Sio$. First, if $t$ is a  finite binary tree labelled in $\Si$, we shall denote by $N_t$ the clopen set $t\cdot T_\Si^\om$. 
Notice that it is easy to see that one can take, as a restricted  basis for the Cantor topology on $T_\Si^\om$, the clopen sets of the form $t_0\cdot T_\Si^\om$, where $t_0$ is a finite labelled binary tree 
whose domain is of the special form $\{l, r\}^{\leq n}$.  

 The Borel hierarchy and the projective hierarchy on $T_\Si^\om$ are defined in the same manner as in the case of the topological space $\Si^\om$. 

 The $\om$-language $\mathbb{P}_\infty\! =\! (0^\star\cdot 1)^\omega$ is a well known example of 
${\bf \Pi}^0_2 $-complete subset of $2^\om$ (see Exercise 23.1 in \cite{Kechris94}). It is the set of 
$\om$-words over $2$ having infinitely many occurrences of the letter $1$. 
Its  complement $2^\om - (0^\star\cdot 1)^\om$ is a 
${\bf \Si}^0_2 $-complete subset of $2^\om$.

 It follows from the definition of the B\"uchi acceptance condition for infinite trees that each tree language recognized by some (non deterministic) B\"uchi tree automaton is an analytic set. 

Niwi\'nski showed that some B\"uchi     recognized   tree languages  are actually  ${\bf \Si}^1_1$-complete sets, \cite{Niwinski85}.  
An example is any tree language $T \subseteq T_\Si^\om$ of the form $\exists \mathrm{Path}(L)$, 
where $L\subseteq \Si^\om$ is a regular $\om$-language which is a ${\bf \Pi}^0_2$-complete subset of $\Si^\om$. 
In particular, for $\Si =2$, the tree language ${\mathcal{L}=\exists \mathrm{Path}(\mathbb{P}_\infty )}$ is 
${\bf \Si}^1_1$-complete and hence non Borel \cite{Niwinski85, PerrinPin,Simonnet92}.

 Notice that its complement   $\mathcal{L}^-=\fa \mathrm{Path}( 2^\om - (0^\star\cdot 1)^\om )$ is a 
${\bf \Pi}^1_1$-complete set. It cannot be accepted by some B\"uchi tree automaton because it is not a ${\bf \Si}^1_1$ set. On the other hand, it can be easily seen that it is accepted by some deterministic Muller tree automaton. 

 We now consider the topology on the space $T_\Si^\om$ generated by the regular languages of trees accepted by some B\"uchi  tree automaton. 

 We prove a version of Theorem \ref{main} as a first step towards the proof that the B\"uchi  topology on a space $T_\Si^\om$ is strong Choquet. We set 
$$\mathbb{T}_\infty := \{ t \in T_2^\om \mid \mbox{ for every path $p$ of } t \;\; \forall k  \geq 0 \;\; \exists i \geq k \;\; p(i) = 1\} .$$
This  set is simply the set of infinite trees over the alphabet $2$  having infinitely many letters $1$ on every (infinite) path. We will work in the spaces of the form  $T_\Si^\om$, where
$\Sigma$ is a finite alphabet with at least two elements. We consider a topology
$\tau_\Sigma$ on  $T_\Si^\om$, and a basis $\mathbb{B}_\Sigma$ for
$\tau_\Sigma$.  We consider the following properties of the family
$(\tau_\Sigma ,\mathbb{B}_\Sigma )_\Sigma$, using the previous
identification:

\begin{enumerate}[label=(P\arabic*)]
\item $\mathbb{B}_\Sigma$ contains the usual basic clopen sets $N_t$,\smallskip

\item $\mathbb{B}_\Sigma$ is closed under finite unions and
  intersections,\smallskip

\item $\mathbb{B}_\Sigma$ is closed under projections, in the sense that
  if $\Gamma$ is a finite set with at least two elements and
  $L \in \mathbb{B}_{\Sigma\times\Gamma}$, then
  $\pi_0[L] \in \mathbb{B}_\Sigma$,\smallskip

\item for each $L \in \mathbb{B}_\Sigma$ there is a closed subset
  $C$ of $T_\Si^\om \times \mathbb{T}_\infty$ (i.e., $C$ is the
  intersection of a closed subset of the Cantor space
  $T_\Si^\om \times T_2^\om$ with
  $T_\Si^\om \times \mathbb{T}_\infty$), which is in
  $\mathbb{B}_{\Sigma\times 2}$, and such that $L =\pi_0[C]$.
\end{enumerate}

 Consider now the set of trees $\mathbb{T}_\infty$. It is easy to see that  $\mathbb{T}_\infty$ is accepted by some   
 {\it deterministic } B\"uchi  tree automaton. On the other hand it is well known that the tree languages accepted by some {\it deterministic} B\"uchi tree automaton are ${\bf \Pi}^0_2$ sets, see \cite{ADMN}. Thus the set 
$\mathbb{T}_\infty$ is actually a ${\bf \Pi}^0_2$ set, it is the intersection of a countable sequence 
$(O_i)_{i\in\mathbb{N}}$ of open sets. We may assume, without loss of generality, that the sequence 
$(O_i)_{i\in\mathbb{N}}$ is decreasing with respect to the inclusion relation. Moreover, each open set $O_i$ is a countable union of basic clopen sets $N_{t_{i,j}}$, $j\geq 0$,  and we may also assume, without loss of generality,  that for all integers $i\geq 0$, and all $j\geq 0$, the finite tree $t_{i,j} \subseteq  \{l, r\}^\star$ has a domain of the form   $\{l, r\}^{\leq n}$ for some integer $n$ greater than $i$. We now state the following result, which is a version of Theorem \ref{main} in the case of trees. 

\begin{thm} \label{main2}  
  Assume that the family $(\tau_\Sigma ,\mathbb{B}_\Sigma )_\Sigma$
  satisfies Properties (P1)-(P4). Then the topologies $\tau_\Sigma$ are
  strong Choquet.
\end{thm}

\noindent\bf Proof.\rm\ 
  We first describe a strategy $\tau$ for Player~2. Player~1 first plays
  $t_0 \in T_\Si^\om$ and a $\tau_\Sigma$-open neighborhood
  $U_0$ of $t_0$. Let $L_0$ in $\mathbb{B}_\Sigma$ with
  $t_0 \in L_0\subseteq U_0$. Property (P4) gives $C_0$ with
  $L_0 =\pi_0[C_0]$. This gives $\alpha_0 \in \mathbb{T}_\infty$ such that
  $(t_0,\alpha_0) \in C_0$. We choose $l^0_0 \in \mathbb{N}$ big
  enough to ensure that if 
  $$s^0_0 := \alpha_0\vert l^0_0\mbox{,}$$ 
  then $ N_{s^0_0}$ is included in the open set $O_1$. We set $w_0 := t_0\vert 1$ and
  $V_0 := \pi_0[C_0\cap (N_{w_0}\times N_{s^0_0})]$. By Properties
  (P1)-(P3), $V_0$ is in $\mathbb{B}_\Sigma$ and thus
  $\tau_\Sigma$-open. Moreover,
  $t_0 \in V_0\subseteq L_0\subseteq U_0$, so that Player~2 respects
  the rules of the game if he plays $V_0$.

  Now Player~1 plays $t_1 \in V_0$ and a $\tau_\Sigma$-open
  neighborhood $U_1$ of $t_1$ contained in $V_0$. Let $L_1$ in
  $\mathbb{B}_\Sigma$ with $t_1 \in L_1\subseteq U_1$. Property (P4)
  gives $C_1$ with $L_1 =\pi_0[C_1]$. This gives
  $\alpha_1 \in \mathbb{T}_\infty$ such that $(t_1,\alpha_1) \in
  C_1$. We choose $l^1_0 \in \mathbb{N}$ big enough to ensure that if
  $s^1_0 := \alpha_1\vert l^1_0$, then $N_{s^1_0}$ is included 
  in the open set $O_1$. 
   As $t_1 \in V_0$, there is
  $\alpha'_0 \in \mathbb{T}_\infty$ such that
  $(t_1,\alpha'_0) \in C_0\cap (N_{w_0}\times N_{s^0_0})$. We choose
  $l^0_1 > l^0_0$ big enough to ensure that if
  $s^0_1 := \alpha'_0\vert l^0_1$, then $s^0_1$ is such that $N_{s^0_1}$ is included in the open set $O_2$. 
   We set $w_1 := t_1\vert 2$ and
  $V_1 := \pi_0[C_0\cap (N_{w_1}\times N_{s^0_1})]\cap\pi_0[C_1\cap
  (N_{w_0}\times N_{s^1_0})]$. Here again, $V_1$ is
  $\tau_\Sigma$-open. Moreover, $t_1 \in V_1\subseteq U_1$ and Player
  2 can play $V_1$.

  Next, Player~1 plays $t_2 \in V_1$ and a $\tau_\Sigma$-open
  neighborhood $U_2$ of $t_2$ contained in $V_1$. Let $L_2$ in
  $\mathbb{B}_\Sigma$ with $t_2 \in L_2\subseteq U_2$. Property (P4)
  gives $C_2$ with $L_2 =\pi_0[C_2]$. This gives
  $\alpha_2 \in \mathbb{T}_\infty$ such that $(t_2,\alpha_2) \in
  C_2$. We choose $l^2_0 \in \mathbb{N}$ big enough to ensure that if
  $s^2_0 := \alpha_2\vert l^2_0$, then the basic open set $N_{s^2_0}$ is included 
  in the open set $O_1$.
   As $t_2 \in V_1$, there is
  $\alpha'_1 \in \mathbb{T}_\infty$ such that
  $(t_2,\alpha'_1) \in C_1\cap (N_{w_0}\times N_{s^1_0})$. We choose
  $l^1_1 > l^1_0$ big enough to ensure that if
  $s^1_1 := \alpha'_1\vert l^1_1$, then the basic open set $N_{s^1_1}$ 
  is included in the open set $O_2$. 
   As $t_2 \in V_1$, there is
  $\alpha''_0 \in \mathbb{P}_\infty$ such that
  $(\sigma_2,\alpha''_0) \in C_0\cap (N_{w_1}\times N_{s^0_1})$. We choose
  $l^0_2 > l^0_1$ big enough to ensure that if
  $s^0_2 := \alpha''_0\vert l^0_2$, then the basic open set $N_{s^0_2}$ 
   is included 
  in the open set $O_3$.
   We set $w_2 := t_2\vert 3$ and
  $V_2 := \pi_0[C_0\cap (N_{w_2}\times N_{s^0_2})]\cap\pi_0[C_1\cap
  (N_{w_1}\times N_{s^1_1})]\cap\pi_0[C_2\cap (N_{w_0}\times
  N_{s^2_0})]$. Here again, $V_2$ is $\tau_\Sigma$-open. Moreover,
  $t_2 \in V_2\subseteq U_2$ and Player~2 can play $V_2$.

  If we go on like this, we build $w_k \in \Si^{\{l, r\}^{\leq k+1}}$ and
  $s^n_l \in 2^{\{l, r\}^\star}$ such that $w_0\subseteq w_1\subseteq ...$ and
  $$s^n_0\!\subsetneqq\! s^n_1\!\subsetneqq\! ...$$ 
  This allows us to define
  $\sigma := \mbox{lim}_{l\rightarrow\infty}~w_l \in T_\Si^\om$
  and, for each $n \in \mathbb{N}$,
  $\beta_n :=\mbox{lim}_{l\rightarrow\infty}~s^n_l \in
  T_2^\om$. Note that $\beta_n \in \mathbb{T}_\infty$ since  
  the basic open set $N_{s^n_l}$ is included 
  in the open set $O_{l+1}$.
  $(\sigma ,\beta_n)$ is the limit of $(w_l,s^n_l)$ as $l$ goes to
  infinity and $N_{w_l}\times N_{s^n_l}$ meets $C_n$ (which is closed in
  $T_\Si^\om \times \mathbb{T}_\infty$),
  $(\sigma ,\beta_n) \in C_n$. Thus
$$\sigma \in \bigcap_{n\in\mathbb{N}}~\pi_0[C_n] =
\bigcap_{n\in\mathbb{N}}~L_n\subseteq \bigcap_{n\in\mathbb{N}}~U_n\subseteq
\bigcap_{n\in\mathbb{N}}~V_n\mbox{,}$$ so that $\tau$ is winning for Player
2.
\hfill{$\square$}\bigskip

We now check that the  B\"uchi  topology on a space $T_\Si^\om$ satisfies Properties (P1)-(P4).

\begin{enumerate}[label=(P\arabic*)]
\item  It is very easy to see that for each finite tree $t$ labelled in $\Si$, there exists a  B\"uchi tree automaton accepting  the usual basic clopen set $N_t$.\smallskip

\item $\mathbb{B}_\Sigma$ is closed under finite unions, because any basic open set in the  B\"uchi topology is accepted by some {\it non-deterministic}  B\"uchi tree automaton. Moreover one can easily show, using a classical product construction, that the class of tree languages accepted by some B\"uchi tree automaton is closed under 
finite  intersections. Thus  $\mathbb{B}_\Sigma$ is closed under finite intersections.\smallskip

\item  It follows easily, from the fact that any basic open set in the  B\"uchi topology is accepted by some 
{\it non-deterministic}  B\"uchi tree automaton, that $\mathbb{B}_\Sigma$ is closed under projections.\smallskip

\item  This property follows from the following lemma, which is very similar to Lemma \ref{lemm} above. 
  
\end{enumerate}

  \begin{lem} \label{lemm2} 
  Let $\Sigma$ be a finite set with at least two elements, and
  $L\subseteq T_\Si^\om$ be a regular tree language accepted by some B\"uchi tree automaton. Then
  there is a closed subset $C$ of
  $T_\Si^\om \times\mathbb{T}_\infty$, which is accepted by some B\"uchi tree automaton
  as a subset of $T_{(\Sigma\times 2)}^\om$ identified with
  $ T_\Si^\om\times  T_2^\om$, and such that $L =\pi_0[C]$.
\end{lem}

\noindent\bf Proof.\rm\ 
  Let $\mathcal{A}=(\Si,Q,q_0,Q_f,\Delta )$   be a B\"uchi tree automaton,
  and $L = L(\mathcal{A})$ be its set of accepted trees.  We call $\chi_f$ the characteristic function of $Q_f$.  It maps the 
  state~$q$ to $1$ if $q \in Q_f$, and to~$0$ otherwise.  The
  function~$\chi_f$ is extended to $T_Q^\om$ by setting
  $t' = \chi_f(t)$ and $t'(s) = \chi_f\big( t(s)\big)$.
  Note that a run~$\rho$ of~$\mathcal{A}$ is accepting if and only if
  $\chi_f(\rho) \in \mathbb{T}_\infty$. Let $C$ be the subset of $T_\Si^\om \times
  \mathbb{T}_\infty$ defined by
  \begin{displaymath}
    C := \big\{ (t, t') \in T_\Si^\om \times
  \mathbb{T}_\infty
         \mid\exists \rho 
         \text{ run of $\mathcal{A}$ on $t$ \text{such that}
                $t' = \chi_F(\rho)$}
         \big\}.
  \end{displaymath}
By definition of $C$, $L =\pi_0[C]$.  Let $K$ be the subset of
  $T_\Si^\om \times T_2^\om \times T_Q^\om$ 
  defined by
  \begin{displaymath}
    K := \big\{ (t,t',\rho) \in T_\Si^\om \times T_2^\om \times T_Q^\om 
         \mid \rho 
         \text{ is a run of $\mathcal{A}$ on $t$ \text{such that} 
                $t' = \chi_F(\rho)$}
         \big\}.
  \end{displaymath}
  As $K$ is compact as a closed subset of the  compact space   $T_\Si^\om \times T_2^\om \times T_Q^\om$           and
  $$C =\pi_{T_\Si^\om \times T_2^\om}[K]\cap (T_\Si^\om \times
  \mathbb{T}_\infty)\mbox{,}$$ 
  the subset~$C$ is a closed subset of
  $T_\Si^\om \times
  \mathbb{T}_\infty$.
  Moreover, it is easy to construct a  B\"uchi tree automaton accepting the tree language $C$. 
\hfill{$\square$}

\begin{cor}
 Let $\Sigma$ be a finite set with at least two elements. Then the B\"uchi topology on $T_\Si^\om$ is strong Choquet. 
\end{cor}

\noindent\bf Proof.\rm\ 
This follows from the fact that the  B\"uchi topology on $T_\Si^\om$ satisfies Properties  (P1)-(P4), and from Theorem \ref{main2}. 
\hfill{$\square$}\bigskip

On the other hand, as in the case of the  B\"uchi topology on $\Si^\om$, the B\"uchi topology on $T_\Si^\om$ is  second countable since there are only countably many possible B\"uchi tree automata (up to  identifications), and it is $T_1$ since it is finer than the usual Cantor topology by Property (P1). However, the B\"uchi topology on 
$T_\Si^\om$ is not Polish, by the following result. 

\begin{thm}
Let $\Sigma$ be a finite set with at least two elements. Then the B\"uchi topology on $T_\Si^\om$ is not metrizable and thus not Polish. 
\end{thm}

\noindent\bf Proof.\rm\ 
Recall that in a metrizable topological space, every closed set is a countable intersection of open sets.  We now show that the  B\"uchi topology on $T_\Si^\om$ does not satisfy this property. We have already recalled that the set $\mathcal{L}=\exists \mathrm{Path}(\mathbb{P}_\infty )$ is ${\bf \Si}^1_1$-complete for the usual topology, and it is open for the  B\"uchi topology since it is accepted by some B\"uchi tree automaton. 
Its complement $\mathcal{L}^-$  is the set of trees in $T_\Si^\om$ having all their paths in $2^\om  \setminus (0^\star\cdot 1)^\om$; it is 
${\bf \Pi}^1_1$-complete  for the usual topology and closed  for the  B\"uchi topology. 

On the other hand every tree language accepted by some B\"uchi tree automaton  is an effective analytic set, i.e., a 
$\Si^1_1$ set, and thus also a (boldface) ${\bf \Si}^1_1$ set (for the usual Cantor topology). 
Moreover, every open set for the  B\"uchi topology on $T_\Si^\om$ is a countable union of basic open sets, and thus a countable union of ${\bf \Si}^1_1$ sets. But the class ${\bf \Si}^1_1$ is closed under countable unions (see \cite{Kechris94}). Therefore, every open set for the B\"uchi topology is a ${\bf \Si}^1_1$ set for the usual topology. 

Towards a contradiction, assume now  that the set $\mathcal{L}^-$ is a  countable intersection of open sets  for the  B\"uchi topology. Then it  is a  countable intersection of  ${\bf \Si}^1_1$ sets for the usual topology. But the class 
${\bf \Si}^1_1$ is closed under  countable intersections and thus  $\mathcal{L}^-$ would be also a ${\bf \Si}^1_1$ set for the usual topology. But $\mathcal{L}^-$  is ${\bf \Pi}^1_1$-complete and thus in 
${\bf \Pi}^1_1 \setminus {\bf \Si}^1_1$ (see \cite{Kechris94}), which is absurd.\hfill{$\square$}
  
\begin{Rem} One can infer,  from the previous results on the B\"uchi topology  on  the space  $T_\Si^\om$ and from Theorem \ref{choquet}, that the B\"uchi topology  on  the space  $T_\Si^\om$ is not regular.\end{Rem}

\begin{Rem} The automatic topology on the space $T_\Si^\om$, which can be defined as in the case of the space $\Si^\om$, 
 is Polish, and the proof of this fact is very similar to the one in the case of the space $\Si^\om$ (see Lemma \ref{lemfer}).\end{Rem}

 We now consider the topology on the space $T_\Si^\om$ generated by the class of all regular languages of trees accepted by some Muller tree automaton. We shall call this topology the Muller topology on $T_\Si^\om$. This topology is clearly $T_1$ since it is finer than the usual topology on $T_\Si^\om$. It is second countable since there are only countably many Muller automata.  It is zero-dimensional,  and thus also regular, 
because the class of regular tree languages over the alphabet $\Si$ is closed under taking complements.  We now recall the following Urysohn metrization theorem (see \cite{Kechris94}). 

\begin{thm} Let $X$ be a second-countable topological space. Then $X$ is metrizable if and only if  its topology is $T_1$  and regular.\end{thm}

 This implies that the Muller topology is  metrizable.  Notice that one can define a distance compatible with this topology in a way very  similar to the way we defined a distance compatible with the B\"uchi topology on $\Sio$. On the other hand we recall the following Becker theorem (see Theorem 4.2.6 in \cite{Gao}). 

\begin{thm} Let $\tau$ be a Polish topology on $X$ and $\tau'$ be a second-countable strong Choquet  topology on $X$  finer than $\tau$. Then every  $\tau'$-open set is ${\bf \Si}_1^1$ in $\tau$.\end{thm}

 This implies that the Muller topology is not strong Choquet because there exist some ${\bf \Pi}_1^1$-complete (for the usual topology), and hence non-${\bf \Si}_1^1$  (for the usual topology),  regular set of trees. Such a regular set of trees is open for the Muller topology, which is finer than the usual topology, but is not ${\bf \Si}_1^1$  for the usual topology. We now summarize the results in this section. 

\begin{thm}  Let $\Si$ be a finite  alphabet having at least two letters. 
\begin{enumerate}
\item   The B\"uchi topology on the space $T_\Si^\om$ is strong Choquet, but it is not regular (and hence not zero-dimensional)  and not metrizable. 

\item  The Muller  topology on the space $T_\Si^\om$   is zero-dimensional,  regular and metrizable,  but it is not strong Choquet. 
\end{enumerate}
\noindent In particular,  the B\"uchi topology and the Muller  topology on $T_\Si^\om$  are not Polish. 
\end{thm}

 If we are just interested in the non-polishness of these topologies, we can argue in a more direct way as follows. We first prove the following proposition. 

\begin{prop}\label{2polish-top}
Let $(X, \tau)$ be a Polish topological space, and let $\tau'$ be another Polish topology on $X$ finer than $\tau$. Then the two topologies $\tau$ and $\tau'$ have the same Borel sets. 
\end{prop}

\noindent\bf Proof.  \rm Let $Id: (X, \sigma) \rightarrow (X, \tau)$ be the identity function on $X$, where the domain $X$ is equipped with the topology $\sigma$ and the range $X$ is equipped with the topology $\tau$. This function is continuous since the topology $\sigma$  is finer than the topology $\tau$.  Notice that this implies, by an easy induction on the rank of a Borel set, that the preimage of any Borel set of $(X, \tau)$ is a Borel set of $(X, \sigma)$, i.e. that every Borel set of  $(X, \tau)$ is a Borel set of $(X, \sigma)$. On the other hand, it follows from Lusin-Suslin's Theorem, see \cite[Theorem 15.1]{Kechris94}, that the (injective) image by the function $Id$  of any Borel set of $(X, \sigma)$ is a Borel set of $(X, \tau)$,  i.e. that every Borel set of  $(X, \sigma)$ is a Borel set of $(X, \tau)$.\hfill{$\square$}

\begin{cor}
Let $\Sigma$ be a finite set with at least two elements. Then the B\"uchi topology and the Muller  topology on $T_\Si^\om$  are not Polish. 
\end{cor}

\noindent\bf Proof.\rm\ 
The B\"uchi topology and the Muller  topology on $T_\Si^\om$ are finer than the usual Cantor topology on $T_\Si^\om$. On the other hand both the B\"uchi topology and the Muller  topology on $T_\Si^\om$ contain some open sets which are ${\bf \Si}_1^1$-complete and hence non Borel (for the usual topology), like,  in the case  $\Si=2$, the set  $\mathcal{L}=\exists \mathrm{Path}( \mathcal{\mathbb{P}_\infty} )$ of infinite trees in $T_\Si^\om$ having at least one  path in the $\om$-language $\mathbb{P}_\infty =(0^\star\cdot 1)^\om$. The conclusion now follows from Proposition \ref{2polish-top}.\hfill{$\square$}

\section{Consequences for our topologies}

\subsection{Consequences not directly related to the polishness, concerning isolated points}

\subsubsection{The case of a space of infinite words}\label{inf-words}

 \hs \noindent \bf Notation.\rm\ If $z \in \{ C,B,A,\delta\}$, then the space 
$(\Sigma^{\mathbb{N}},\tau_z)$ is denoted $\mathcal{S}_z$. The set of ultimately periodic $\omega$-words on $\Sigma$ is denoted 
$\mbox{Ult} := \big\{ u\cdot v^\omega \mid u,v \in \Sigma^*\!\setminus\!\{\emptyset\}\big\}$, and 
$P := \Sigma^{\mathbb{N}}\!\setminus\!\mbox{Ult}$.

\noindent (1) As noted in \cite{SchwarzS10}, Ult is the set of \emph{isolated points} of $\mathcal{S}_B$ and $\mathcal{S}_A$ (recall that a point 
$\sigma\in\Sigma^{\mathbb{N}}$ is isolated if $\{\sigma\}$ is an open set). Indeed, each singleton 
$\{ u\cdot v^\omega\}$\! formed by an ultimately periodic $\omega$-word is an $\omega$-regular language, and thus each ultimately periodic $\omega$-word is a $\mathcal{S}_A$-isolated point. Conversely, if $\{\sigma\}$ is $\tau_B$-open, then it is $\omega$-regular and then the $\omega$-word $\sigma$ is ultimately periodic (because any countable $\omega$-regular language contains only ultimately periodic $\omega$-words, see \cite{Buchi62,PerrinPin,Staiger97}).

\noindent (2) Every nonempty $\omega$-regular set contains an ultimately
periodic $\omega$-word, \cite{Buchi62,PerrinPin,Staiger97}. In particular, the set Ult of
isolated points of $\mathcal{S}_B$ and $\mathcal{S}_A$ is \emph{dense},
and a subset of $\mathcal{S}_B$ or $\mathcal{S}_A$ is dense if and only if
it contains Ult.

\subsubsection{The case of a space of infinite trees}\indent

Recall that a tree  $t\in T_\Si^\om$    is    regular  if and only if  for each $a\in \Si$ the set $\{u \in \{l, r\}^\star \mid t(u)=a \}$ is a regular set of finite words over the alphabet $\{l, r\}$. 
For each regular tree $t\in T_\Si^\om$, the singleton $\{t\}$ is a (closed) regular tree language accepted by some B\"uchi tree automaton. 
Moreover, a regular tree language accepted by some Muller or Rabin tree automaton is non-empty if and only if it contains a regular tree, see \cite{Thomas90}.  

Therefore we can state properties similar to those stated in the case of a space of words in Section \ref{inf-words}.  For the automatic, the B\"uchi or the Muller topologies on a space $ T_\Si^\om$, the set of isolated points is the set $Reg$-$trees$ of regular trees, and this set is dense.  It is left here to the reader to see how Properties $(1)$-$(2)$ of the preceding section are extended to the case of a space of trees (in a very similar way). 

\subsection{Consequences of the polishness}\indent

 Here we concentrate on the case of a space of infinite words. We consider our topologies on 
$\Sigma^{\mathbb{N}}$, where $\Sigma$ is a finite set with at least two elements. We refer to \cite{Kechris94} when classical descriptive set theory is involved. 

\noindent (a) The union $P\cup\mbox{Ult}$ is the \emph{Cantor-Bendixson
  decomposition} of $\mathcal{S}_B$ and $\mathcal{S}_A$ (see Theorem 6.4 in
\cite{Kechris94}). This means that $P$ is perfect (i.e., closed without
isolated points) and Ult is countable open. Let us check that $P$ is
perfect. We argue by contradiction, so that we can find $\sigma \in P$ and
an $\omega$-regular language $L$ such that
$\{\sigma\} = L\!\setminus\!\mbox{Ult}$. Note that
$L\subseteq \{\sigma\}\cup\mbox{Ult}$ is countable. But a countable regular
$\om$-language contains only ultimately periodic words (see
\cite{PerrinPin}), and thus $L\subseteq \mbox{Ult}$, which is
absurd.

\noindent (b) The closed subspace $(P,\tau_B)$ of $\mathcal{S}_B$ is
\emph{homeomorphic to the Baire space} $\mathbb{N}^\mathbb{N}$. Indeed, it
is not empty since Ult is countable and $\Sigma^\mathbb{N}$ is not,
zero-dimensional and Polish as a closed subspace of the zero-dimensional
Polish space $\mathcal{S}_B$. By Theorem 7.7 in \cite{Kechris94}, it is enough to
prove that every compact subset of $(P,\tau_B)$ has empty interior. We
argue by contradiction, which gives a compact set $K$. Note that there is
an $\omega$-regular language $L$ such that $P\cap L$ is a nonempty compact
subset of $K$, so that we may assume that $K = P\cap L$. Theorem \ref{buch}
gives $(U_i)_{i<n}$ and $(V_i)_{i<n}$ with
$L =\bigcup_{i<n}~U_i\cdot V_i^\omega$.  On the other hand, $L$ is not
countable since every countable regular $\om$-language contains only
ultimately periodic words and $K = P\cap L$ is non-empty. Thus $n > 0$
and, for example, $U_0\cdot V_0^\omega$ is not countable.

This implies that we can find $v_0,v_1 \in V_0$ which are not powers of the
same word (which means that we cannot find a finite word $u$ and integers $m,n$ with $v_0=u^m$ and $v_1=u^n$).  Indeed, we argue by contradiction. Let
$v_0 \in V_0 \setminus \{\emptyset\}$, and $v \in \Sigma^*$ of minimal
length such that $v_0$ is a power of $v$. We can find
$w_0,w_1,... \in V_0\!\setminus\!\{\emptyset\}$ such that
$\sigma := w_0\cdot w_1... \neq  v^\om$.  Fix a natural number $i$. Then $w_i$
and $v_0$ are powers of the same word $w$. By Corollary 6.2.5 in
\cite{Lothaire02}, $v$ and $w$ are powers of the same word $u$, and $v_0$
too. By minimality, $u = v$, and $w_i$ is a power of $v$. Thus
$\sigma = v^\om$, which is absurd.

Let $u_0 \in U_0$, and $L' := \{ u_0\}\cdot\{ v_0,v_1\}^\omega$. Note
that $P\cap L'$ is a $\tau_B$-closed subset of $K$, so that it is
$\tau_B$-compact. As the identity map from $(P\cap L',\tau_B)$ onto $(P\cap L',\tau_C)$ is
continuous, $P\cap L'$ is also $\tau_C$-compact. But the map
$\alpha\!\mapsto\! u_0\cdot v_{\alpha (0)}\cdot v_{\alpha (1)}\ldots$ is a homeomorphism
from the Cantor space onto $L'$, by Corollaries 6.2.5 and 6.2.6 in
\cite{Lothaire02}. Thus $P\cap L' = L'\!\setminus\!\mbox{Ult}$ is a dense
closed subset of $L'$. Thus $P\cap L' = L'$, which is absurd since
$u_0\cdot v_0^\om \in L'\!\setminus\! P$.

\section{Concluding remarks}\indent

We obtained in this paper new links and interactions between descriptive set theory and theoretical computer science, showing that two topologies considered in \cite{SchwarzS10} are Polish.

 Notice that this paper is also motivated by the fact that the Gandy-Harrington topology, generated by the effective analytic subsets of a recursively presented Polish space, is an extremely powerful tool in descriptive set theory. In particular, this topology is
  used to prove some results of classical type (without reference to effective descriptive set theory in their statement). Among these results, let us mention the dichotomy theorems in \cite{HKL,KST99,Lec09a,Lec13}. Sometimes, no other proof is known. Part of the power of this technique comes from the nice closure properties of the class $\Ana$ of effective analytic sets (in particular the closure under projections). 
 
  The class of $\omega$-regular languages has even stronger closure properties. So our hope is that the study of the B\"uchi topology, generated by the $\omega$-regular languages, will help to prove some automatic versions of known descriptive results in the context of theoretical computer science. For instance, more precisely, let $\Sigma,\Gamma$ be finite sets with at least two elements, and $L$ be a subset of $\Sigma^{\mathbb{N}}\!\times\!\Gamma^{\mathbb{N}}$ which is $\omega$-regular and also a countable union of Borel rectangles. It would be very interesting to know whether $L$ is open for the product topology $\tau_B\!\times\!\tau_B$. Indeed, this would give a version of the ${\mathcal G}_0$-dichotomy for $\om$-regular languages, and thus a very serious hope to get versions of many difficult dichotomy results of descriptive set theory for $\om$-regular languages (see \cite{Miller12}).

 From Theorem~\ref{thm:main}, we know that there is a complete distance which is compatible with $\tau_z$. It would be interesting to have a natural complete distance compatible with $\tau_B$. We leave this as an open question for further study.

  \hs {\bf Acknowledgements.}
We are  very grateful to the anonymous referees for their careful reading and their useful comments on a preliminary version of our paper that led to a great improvement of the presentation of the paper. The second author thanks very much Henryk Michalewski who suggested to study the analogue of the  B\"uchi topology in the case of a space of infinite labelled trees.

\end{document}